\documentclass[11pt,reqno]{amsart}
\usepackage{amssymb,amsthm,amsmath,verbatim,enumerate,ifthen}
\usepackage[mathscr]{eucal}
\usepackage[utf8]{inputenc}
\usepackage[T1]{fontenc}
\def\N{\mathbb{N}}
\def\R{\mathbb{R}}

\def\QQ{\mathbb{Q}}

\def\Z{\mathbb{Z}}

\def\C{\mathscr{C}}

\def\E{\mathscr{E}}

\def\T{\mathscr{T}}

\def\cl{\mathop{\mbox{\rm cl}}}
\def\conv{\mathop{\mbox{\rm conv}}\nolimits}

\def\div{\mathop{\text{\rm div}}\nolimits}
\def\diam{\mathop{\text{\rm diam}}\nolimits}
\def\Cdot{\!\cdot\!}

\newtheorem{theorem}{Theorem}[section]
\newtheorem*{theorem*}{Theorem}
\def\Thm#1#2{\ifthenelse{\equal{#1}{*}}{\begin{theorem*}#2\end{theorem*}}
             {\begin{theorem}\label{T#1}#2\end{theorem}}}
\newtheorem{Atheorem}{Theorem}

\def\thm#1{Theorem~\ref{T#1}}

\newtheorem{proposition}[theorem]{Proposition}
\newtheorem*{proposition*}{Proposition}
\def\Prp#1#2{\ifthenelse{\equal{#1}{*}}{\begin{proposition*}#2\end{proposition*}}
             {\begin{proposition}\label{P#1}#2\end{proposition}}}

\newtheorem{corollary}[theorem]{Corollary}
\newtheorem*{corollary*}{Corollary}
\def\Cor#1#2{\ifthenelse{\equal{#1}{*}}{\begin{corollary*}#2\end{corollary*}}
             {\begin{corollary}\label{C#1}#2\end{corollary}}}
\def\cor#1{Corollary~\ref{C#1}}

\newtheorem{lemma}[theorem]{Lemma}
\newtheorem*{lemma*}{Lemma}
\def\Lem#1#2{\ifthenelse{\equal{#1}{*}}{\begin{lemma*}#2\end{lemma*}}
             {\begin{lemma}\label{L#1}#2\end{lemma}}}
\def\lem#1{Lemma~\ref{L#1}}

\theoremstyle{definition}
\newtheorem{remark}[theorem]{Remark}
\newtheorem*{remark*}{Remark}
\def\Rem#1#2{\ifthenelse{\equal{#1}{*}}{\begin{remark*}\rm #2\end{remark*}}
             {\begin{remark}\label{R#1}\rm #2\end{remark}}}

\newtheorem{example}[theorem]{Example}
\newtheorem*{example*}{Example}
\def\Exa#1#2{\ifthenelse{\equal{#1}{*}}{\begin{example*}\rm #2\end{example*}}
             {\begin{example}\label{Ex#1}\rm #2\end{example}}}

\def\eq#1{{\rm(\ref{E#1})}}
\def\Eq#1#2{\ifthenelse{\equal{#1}{*}}
  {\begin{equation*}\begin{aligned}[]#2\end{aligned}\end{equation*}}
  {\begin{equation}\begin{aligned}[]\label{E#1}#2\end{aligned}\end{equation}}}

\begin{document}
\vspace{5mm}

\date{\today}

\title[]{Convexity of sets in metric Abelian groups}

\author[W. Fechner]{W{\l}odzimierz Fechner}
\address{Institute of Mathematics, Lodz University of Technology, ul. W\'olcza\'nska 215, 90-924 \L\'od\'z, Poland}
\email{wlodzimierz.fechner@p.lodz.pl}

\author[Zs. P\'ales]{Zsolt P\'ales}
\address{Institute of Mathematics, University of Debrecen, 
H-4002 Debrecen, Pf.\ 400, Hungary}
\email{pales@science.unideb.hu}


\thanks{The research of the second author was supported by the EFOP-3.6.1-16-2016-00022 and the EFOP-3.6.2-16-2017-00015 projects. These projects are co-financed by the European Union and the European Social Fund.}

\begin{abstract}
In the present paper, we introduce a new concept of convexity which is generated by a family of endomorphisms of an Abelian group. In Abelian groups equipped with a translation invariant metric, we define the boundedness, the norm, the measure of injectivity and the spectral radius of endomorphisms. Beyond the investigation of their properties, our first main goal is an extension of the celebrated R\aa dstr\"om Cancellation Theorem. Another result generalizes the Neumann Invertibility Theorem. Next we define the convexity of sets with respect to a family of endomorphisms and we describe the set-theoretical and algebraic structure of the class of such sets. Given a subset, we also consider the family of endomorphisms that make this subset convex and we establish the basic properties of this family. Our first main result establishes conditions which imply midpoint convexity. The next main result, using our extension of the R\aa dstr\"om Cancellation Theorem, presents further structural properties of the family of endomorphisms that make a subset convex.
\end{abstract}

\maketitle

\section{Introduction}

The notion of convexity is usually defined in the context of linear spaces. If $X$ is a linear space, then a subset $D\subseteq X$ is termed \emph{convex} if, for all $x,y\in D$ and $t\in[0,1]$,
\Eq{tD}{
  tx+(1-t)y\in D.
}
In other words, $D$ is convex if it is closed with respect to two-term convex combinations. This notion has been generalized in several ways,  for instance, given a fixed number $t\in[0,1]$, a set $D$ is called \emph{$t$-convex} if the inclusion \eq{tD} is valid for all $x,y\in D$. In particular, the $\frac12$-convex sets are said to be \emph{midpoint convex}. It is clear that these classes of subsets are closed with respect to intersection and therefore, one can naturally define also the convex, $t$-convex and midpoint convex hulls of sets, this is the subject of several papers, e.g. \cite{Rei76,ZhaShi09}. The structural properties of the set of numbers $t\in[0,1]$ such that $D$ is $t$-convex were intensively discussed, see \cite{Kuc85} and some references therein. The notion of convexity in metric spaces was introduced by Menger \cite{Men28}, see also \cite{LelNit60}, approximate convexity properties of sets were investigated in \cite{DilHowRob99,DilHowRob00,DilHowRob02}. The general theory of convex structures was developed in the book \cite{Vel93}. A possible approach to introduce convexity in topological, particularly in Lie groups is given by R\aa dstr\"om \cite{Rad52a}.

Two notions of convexity (which are generalizations of $\QQ$-convexity) have been introduced in the paper \cite{JarPal15} in the Abelian semigroup setting. There, the main motivation was to extend the Stone Decomposition Theorem, which is the key result to prove the classical separation theorems for convex sets in the linear space setting (cf. \cite{Hol72}, \cite{Sto46}). For the separation of subsemigroups, Stone type theorems were obtained by the second author in \cite{Pal89d} and were applied for the characterization of quasideviation means in \cite{Pal89b}. 

In the present paper, we introduce a new concept of convexity which is generated by a family of endomorphisms of an Abelian group. In the case when the underlying group is the additive group of a linear space and the endomorphisms are given as multiplication by scalars from the interval $[0,1]$, these notions reduce to the above mentioned concepts of convexity. The most convenient setting where our results will be formulated are Abelian groups equipped with a translation invariant metric. In this framework, we introduce the boundedness, the norm, the measure of injectivity and the spectral radius of endomorphisms motivated by the respective concepts from the theory of linear operators on normed spaces. Beyond the investigation of their properties, the first main goal of Section 2 is \thm{RCT}, which extends R\aa dstr\"om's celebrated cancellation theorem. Another useful result of this section is \thm{NIT} which generalizes the Neumann Invertibility Theorem. In Section 3, we define the convexity of sets with respect to a family of endomorphisms, then we describe the set-theoretical and algebraic structure of the class of such sets. Given a subset, we also consider the family of endomorphisms that make this subset convex and we establish the basic properties of this family. One of the main results of this section is \thm{2} which establishes conditions which imply the midpoint convexity. The second main result is \thm{nk} which, under certain regularity assumptions and using our extension of the R\aa dstr\"om Cancellation Theorem, presents further structural properties of the family of endomorphisms that make a subset convex.

\section{Auxiliary concepts and tools}

Let $(X,+)$ be an Abelian group and let $\E(X)$ denote the family of all endomorphisms, i.e., additive maps $T:X\to X$. Then, one can see that $\E(X)$ is a ring if the composition is taken as the ring multiplication. Thus, every element $T\in\E(X)$ generates an endomorphism $\widetilde{T}:\E(X)\to\E(X)$ defined by $\widetilde{T}(S):=T\circ S$. For a family $\T\subseteq \E(X)$ we denote by $\widetilde{\T}$ the set $\{\widetilde{T} \mid T \in \T\}$. Finally, $I$ is the identity map of $X$. 

The multiplication of the elements of $X$ by natural numbers can be introduced via the following recursive definition: 
\Eq{*}{
  1\Cdot x:=x,\qquad\mbox{and}\qquad (n+1)\Cdot x:=n\Cdot x+x \qquad(x\in X,\,n\in\N).
}
The mapping $\pi_n(x):=n\Cdot x$ is always an endomorphism provided that $(X,+)$ is an Abelian group. We say that an Abelian group $(X,+)$ is \emph{divisible by} $n\in\N$ if the map $\pi_n$ is a bijection (and hence an automorphism) of $X$. In this case, for $x\in X$, the element $\pi_n^{-1}(x)$ will simply be denoted by $\frac{1}{n}\Cdot x$. The set of natural numbers $n$ for which $X$ is uniquely divisible by $n$ will be denoted by $\div(X)$. Obviously, $\div(X)$ is a multiplicative subsemigroup of $\N$ whose unit element is $1$. 

For the sake of convenience, we introduce the following notation: For a subset $A\subseteq X$ and $n\in\N$, 
\Eq{*}{
  n\Cdot A:=\{n\Cdot x\mid x\in A\}\qquad\mbox{and}\qquad
  [n]A:=\{x_1+\dots+x_n\mid x_1,\dots,x_n\in A\}.
}
The inclusion $n\Cdot A\subseteq [n]A$ is obvious. In the case when the reversed inclusion is also valid we say that the set $A$ is \emph{$n$-convex}. For properties of $n$-convex sets, we refer to the paper \cite{JarPal15}. In particular, by \cite[Proposition 2]{JarPal15}, we have that if a set is $n$- and also $m$-convex, then it is $(nm)$-convex.

In the case when the Abelian group $(X,+)$ is equipped with a \emph{translation invariant metric $d$}, we say that $(X,+,d)$ is a \emph{metric Abelian group}. In a metric Abelian group, the group operations are continuous with respect to the topology induced by $d$. Therefore, metric groups are automatically topological groups. In a metric Abelian group $(X,+,d)$, the \emph{$d$-norm} $\|\Cdot\|_d:X\to\R$ is defined as $\|x\|_d:=d(x,0)$. From the properties of $d$, one can deduce that $\|\Cdot\|_d$ is a positive definite and subadditive even function on $X$ and, for all $x,y\in X$, $d(x,y)=\|x-y\|_d$. The positive homogeneity of the $d$-norm is not a consequence of the properties of a translation invariant metric. In general, the subadditivity of $\|\Cdot\|_d$ implies that $\|n\Cdot x\|_d\leq n\|x\|_d$ for all $x\in X$ and $n\in\N$. The equality here, however, may not be valid.  

An endomorphism $T:X\to X$ is called \emph{$d$-bounded} if there exists $c\geq0$ such that $\|T(x)\|_d\leq c\|x\|_d$ for all $x\in X$. The smallest number $c$ satisfying this condition is called the \emph{$d$-norm} of $T$ and is denoted by $\|T\|^*_d$. The symbol $\E^d(X)$ will denote the subring of $\E(X)$ of all $d$-bounded endomorphisms. More generally, for $\T\subseteq\E(X)$, the symbol $\T^d$ denotes the $d$-bounded elements of $\T$. One can show that $(\E^d(X),+,d^*)$ is also a metric group, where $d^*(T,S):=\|T-S\|^*_d$ for all $T,S\in\E^d(X)$, furthermore, $\|\Cdot\|^*_d$ is submultiplicative, which makes $\E^d(X)$ a $d$-normed ring. The smallest number $c$ such that 
$\|n\Cdot x\|_d\leq c\|x\|_d$ for all $x\in X$, that is $\|\pi_n\|^*_d$, will simply be denoted by $\|n\|^*_d$.

For a $d$-bounded endomorphism $T$, we define its \emph{measure of injectivity $\mu_d(T)$} by
\Eq{*}{
  \mu_d(T):=\inf_{x\in X\setminus\{0\}}\frac{\|T(x)\|_d}{\|x\|_d}.
}
It immediately follows from the $d$-boundedness of $T$ that $\mu_d(T)\leq\|T\|_d^*$. The definition also implies 
\Eq{muT1}{
   \mu_d(T)\|x\|_d\leq \|T(x)\|_d \qquad(T\in\E^d(X),\,x\in X).
}
Therefore, if $\mu_d(T)>0$, then $T$ must be injective. Furthermore, in this case, with the substitution $y=T(x)$, for the inverse map $T^{-1}:T(X)\to X$, we can obtain the following inequality:
\Eq{*}{
  \|T^{-1}(y)\|_d\leq \frac{1}{\mu_d(T)}\|y\|_d \qquad(y\in T(X)).
}
If, additionally, $T$ is surjective, then it follows that $T^{-1}\in\E^d(X)$ and $\|T^{-1}\|^*_d=(\mu_d(T))^{-1}$. Conversely, if $T\in\E^d(X)$ is a bijective map such that $T^{-1}$ is bounded, then $\mu_d(T)=(\|T^{-1}\|^*_d)^{-1}>0$. If $X$ is the additive group of a Banach space, then the open mapping theorem implies that the inverse of bijective linear map is automatically bounded, therefore, its modulus of injectivity is positive.

For $n\in\N$, the measure of injectivity of the map $\pi_n$ will simply be denoted by $\mu_d(n)$. In this particular case we also have
\Eq{mun}{
   \mu_d(n)\|x\|_d\leq \|n\Cdot x\|_d \qquad(n\in\N,\,x\in X).
}

The function $\mu_d:\E^d(X)\to[0,\infty[\,$ enjoys a \emph{supermultiplicativity} and a \emph{Lipschitzian} property as stated below.

\Lem{mu}{Let $(X,+,d)$ be a metric group. Then, for all $T,S\in\E^d(X)$,
\Eq{muT2}{
 \mu_d(T)\|S\|^*_d\leq\|T\circ S\|^*_d,\qquad
  \mu_d(T)\mu_d(S)\leq\mu_d(T\circ S)\qquad\mbox{and}\qquad 
  |\mu_d(T)-\mu_d(S)|\leq\|T-S\|^*_d.
}
In particular, for all $m,n\in\N$,
\Eq{mu}{
  \mu_d(n)\|m\|^*_d\leq\|nm\|^*_d,\qquad
  \mu_d(n)\mu_d(m)\leq\mu_d(nm)\qquad\mbox{and}\qquad 
  |\mu_d(n)-\mu_d(m)|\leq |n-m|.
}}

\begin{proof}
Let $T,S\in\E^d(X)$ be arbitrary. Then, for $x\in X\setminus\{0\}$, the inequality \eq{muT1} implies that
\Eq{*}{
  \mu_d(T)\frac{\|S(x)\|_d}{\|x\|_d}
  \leq \frac{\|(T\circ S)(x)\|_d}{\|x\|_d}.
}
Taking the supremum and the infimum of both sides with respect to $x\in X\setminus\{0\}$, respectively, the two first two inequalities in \eq{muT2} follow.

For the last inequality, observe that, for all $x\in X\setminus\{0\}$,
\Eq{*}{
  \mu_d(T)\leq \frac{\|T(x)\|_d}{\|x\|_d}
  \leq \frac{\|S(x)\|_d+\|T(x)-S(x)\|_d}{\|x\|_d}
  \leq \frac{\|S(x)\|_d}{\|x\|_d}+\|T-S\|^*_d.
}
Taking the infimum of the right hand side, it follows that
\Eq{*}{
 \mu_d(T)\leq \mu_d(S)+\|T-S\|^*_d.
}
Now, interchanging the roles of $S$ and $T$, we obtain the last inequality in \eq{muT2}.

The inequalities in \eq{mu} are immediate consequences of \eq{muT2}.
\end{proof}

The above lemma immediately yields the following result. Its simple proof is left to the reader.

\Cor{mu}{Let $(X,+,d)$ be a metric group. Then the subset
\Eq{*}{
  \E^+_d(X):=\{T\in\E^d(X)\mid \mu_d(T)>0\}
}
is an open multiplicative subsemigroup of $\E^d(X)$.}

Further properties of $d$-bounded maps are stated in the following lemma.

\Lem{nx}{Let $(X,+,d)$ be a metric group and $T\in\E^d(X)$. Then $T$ maps  bounded, resp.\ compact subsets of $X$ into bounded, resp.\ compact subsets. If $\mu_d(T)>0$ and either $X$ is complete or $T(X)$ is closed, then $T$ maps closed subsets of $X$ into closed subsets.}

\begin{proof} If $T\in\E^d(X)$ then $T$ is Lipschitzian and hence continuous (with Lipschitz modulus $\|T\|^*_d$), therefore the first statement follows. 

Now assume that $\mu_d(T)>0$, and $D\subseteq X$ is closed. To show that $T(D)$ is closed, let $(b_k)$ be a convergent sequence from $T(D)$ with a limit point $b_0$. Then there exists a sequence $(d_k)$ in $D$ such that $T(d_k)=b_k$ for all $k\in\N$.

First assume that $X$ is complete. By inequality \eq{muT1}, for all $k,m\in\N$, it follows that
\Eq{*}{
  \mu_d(T)\|d_k-d_m\|_d\leq\|T(d_k)-T(d_m)\|_d=\|b_k-b_m\|.
}
Since $(b_k)$ is a Cauchy sequence, therefore, $(d_k)$ is a Cauchy sequence, too. By the completeness of $X$, the sequence $(d_k)$ must possess a limit $d_0\in D$. The continuity of $T$ yields $T(d_0)=b_0$, which proves that $b_0\in T(D)$.

Secondly, assume that $T(X)$ is closed. Then, for all $k\in\N$, we have $b_k\in T(X)$, which implies that $b_0\in T(X)$. Thus, there exists $d_0\in X$  such that $b_0=T(d_0)$. By inequality \eq{muT1}, for all $k\in\N$, it follows that
\Eq{*}{
  \mu_d(T)\|d_k-d_0\|_d\leq\|T(d_k)-T(d_0)\|_d=\|b_k-b_0\|.
}
Therefore, $(d_k)$ converges to $d_0\in D$, which again proves that $b_0\in T(D)$. 
\end{proof}

The following result is an extension of the celebrated R{\aa}dstr\"{o}m Cancellation Theorem (cf.\cite{Rad52b}).

\Thm{RCT}{Let $(X,+,d)$ be a metric Abelian group and let $n_0\in\N$ such that $\mu_d(n_0)>1$. Let $A\subseteq X$ be an arbitrary subset, let $B\subseteq X$ be closed and $n_0$-convex subset, and $C\subseteq X$ be a $d$-bounded nonempty subset such that $A+C\subseteq B+C$. Then $A\subseteq B$.}

\begin{proof} Assume that $A,B,C$ satisfy the conditions of the theorem. First, we will prove, for all $n\in\N$, by induction that 
\Eq{*}{
  [n]A+C\subseteq [n]B+C.
}
For $n=1$, this is just the assumption $A+C\subseteq B+C$. If it holds for $n$, then
\Eq{*}{
  [n+1]A+C=A+[n]A+C\subseteq A+[n]B+C=[n]B+A+C\subseteq [n]B+B+C=[n+1]B+C.
}
In particular, in view of the $n_0$-convexity of $B$, for $n=n_0^k$, we get
\Eq{*}{
   n_0^k\Cdot A+C\subseteq [n_0^k]A+C\subseteq [n_0^k]B+C=n_0^k\Cdot B+C.
}
Let $c_0\in C$ be fixed and let $a\in A$ be arbitrary. Then, by the above inclusion, for all $k\in\N$, there exist $b_k\in B$ and $c_k\in C$ such that
\Eq{*}{
  n_0^{k}\Cdot a+c_0=n_0^{k}\Cdot b_k+c_k.
}
Hence, applying inequality \eq{mun} with the substitutions $n:=n_0^k$ and $x:=a-b_k$ together with inequality \eq{mu} of \lem{mu}, we get
\Eq{*}{
  +\infty>\diam_d(C)\geq\|c_k-c_0\|_d
  =\|n_0^{k}\Cdot (a-b_k)\|_d
  \geq \mu_d(n_0^k)\|a-b_k\|_d
  \geq \mu_d(n_0)^k\|a-b_k\|_d.
}
By condition $\mu_d(n_0)>1$, it follows that the sequence $\|a-b_k\|_d$ converges to zero. Hence, the sequence $(b_k)$ converges to $a$ and the closedness of $B$ now yields that $a$ must belong to $B$. This proves the desired inclusion $A\subseteq B$.
\end{proof}

The \emph{$d$-spectral radius of an endomorphism $T\in\E^d(X)$} is defined as
\Eq{*}{
  \rho_d(T):=\limsup_{m\to\infty}\sqrt[m]{\|T^m\|^*_d}.
}
In view of the submultiplicativity of $\|\Cdot\|^*_d$ on $\E^d(X)$, for all $m\in\N$, we have that 
\Eq{*}{
  \sqrt[m]{\|T^m\|^*_d}
  \leq \sqrt[m]{(\|T\|^*_d)^m}
  =\|T\|^*_d.
}
Now, upon taking the upper limit as $m\to\infty$, we get
\Eq{*}{
  \rho_d(T)=\limsup_{m\to\infty}\sqrt[m]{\|T^m\|^*_d}
  \leq\|T\|^*_d.
}

The next statement summarizes the basic properties of the spectral radius and its connection to the injectivity modulus.

\Lem{SR}{Let $(X,+,d)$ be a metric group and let $T,S\in\E^d(X)$ such that $T\circ S=S\circ T$. Then
\Eq{SR}{
  \rho_d(T+S)\leq\rho_d(T)+\rho_d(S) \qquad\mbox{and}\qquad
  \mu_d(T)\rho_d(S)\leq \rho_d(T\circ S)\leq \rho_d(T)\rho_d(S).
}
In particular,
\Eq{*}{
  \mu_d\leq\rho_d\leq\|\cdot\|^*_d.
}}

\begin{proof} Assume that $T$ and $S$ are commuting endomorphisms. Then, by the binomial theorem, for all $m\in\N$, we have that
\Eq{TS}{
  (T+S)^m=\sum_{k=0}^m\binom{m}{k}T^k\circ S^{m-k}.
}
Let $\varepsilon>0$ be arbitrary. Then, there exists $n_0\in\N$ such that, for all $n\geq n_0$,
\Eq{TS1}{
   \|T^n\|^*_d<(\rho_d(T)+\tfrac\varepsilon2)^n \qquad\mbox{and}\qquad
   \|S^n\|^*_d<(\rho_d(S)+\tfrac\varepsilon2)^n.
}
Denote
\Eq{*}{
  A(T,S)
  :=\sum_{k=0}^{n_0-1}\binom{m}{k}\Big|\|T^k\|-(\rho_d(T)+\tfrac\varepsilon2)^k\Big|(\rho_d(S)+\tfrac\varepsilon2)^{-k}.
}
With this notation, for $m\geq 2n_0-1$, from \eq{TS} and \eq{TS1} we can deduce
\Eq{*}{
  \|(T+S)^m\|^*_d
  &\leq\sum_{k=0}^m\binom{m}{k}\|T^k\|^*_d\|S^{m-k}\|^*_d \\
  &=\sum_{k=0}^{n_0-1}\binom{m}{k}\|T^k\|^*_d\|S^{m-k}\|^*_d
  +\sum_{k=0}^{n_0-1}\binom{m}{k}\|T^{m-k}\|^*_d\|S^k\|^*_d
  +\sum_{k=n_0}^{m-n_0}\binom{m}{k}\|T^k\|^*_d\|S^{m-k}\|^*_d\\
  &\leq\sum_{k=0}^{n_0-1}\binom{m}{k}
        \|T^k\|^*_d(\rho_d(S)+\tfrac\varepsilon2)^{m-k}
  +\sum_{k=0}^{n_0-1}\binom{m}{k}
        (\rho_d(T)+\tfrac\varepsilon2)^{m-k}\|S^k\|^*_d\\
  &\qquad+\sum_{k=n_0}^{m-n_0}\binom{m}{k}
  (\rho_d(T)+\tfrac\varepsilon2)^k(\rho_d(S)+\tfrac\varepsilon2)^{m-k}\\
  &\leq A(T,S)(\rho_d(S)+\tfrac\varepsilon2)^{m}
  +A(S,T)(\rho_d(T)+\tfrac\varepsilon2)^{m}
  +\sum_{k=0}^{m}\binom{m}{k}
  (\rho_d(T)+\tfrac\varepsilon2)^k(\rho_d(S)+\tfrac\varepsilon2)^{m-k}\\
  &= A(T,S)(\rho_d(S)+\tfrac\varepsilon2)^{m}
  +A(S,T)(\rho_d(T)+\tfrac\varepsilon2)^{m}
  +(\rho_d(T)+\rho_d(S)+\varepsilon)^{m}.
}

Taking the $m$th root of both sides and then (using elementary tools) calculating the upper limit as $m\to\infty$, we get
\Eq{*}{
  \rho_d(T+S)\leq \rho_d(T)+\rho_d(S)+\varepsilon.
}
The number $\varepsilon>0$ being arbitrary, the first inequality in \eq{SR} follows.

To prove the first part of the second inequality in \eq{SR}, using \lem{mu}, observe that
\Eq{*}{
  \mu_d(T)^m\|S^m\|^*_d\leq \|T^m\circ S^m\|^*_d\leq \|(T\circ S)^m\|^*_d.
}
Taking the $m$-th root and calculating the upper limit as $m\to\infty$, we get $\mu_d(T)\rho_d(S)\leq\rho_d(T\circ S)$. The other inequality in \eq{SR} is the consequence of the inequalities
\Eq{*}{
  \|(T\circ S)^m\|^*_d\leq\|T^m\circ S^m\|^*_d
  \leq\|T^m\|^*_d\|S^m\|^*_d.
}
Finally, putting $S:=I$ into \eq{SR}, we get the last assertion.
\end{proof}

The following result is an analogue of the so-called Neumann invertibility theorem and it can be proved almost exactly in the same way as the classical result of Carl Neumann.

\Thm{NIT}{Let $(X,+,d)$ be a complete metric Abelian group and let $T\in\E^d(X)$ such that $\rho_d(T)<1$. Then $I-T$ is an invertible element of $\E^d(X)$, furthermore,
\Eq{*}{
  (I-T)^{-1}=\sum_{k=0}^\infty T^k.
}}

\begin{proof} It easily follows from the completeness of the metric space $(X,d)$ that $\E^d(X)$ is also complete in the metric $d^*$ defined above.

From the condition $\rho_d(T)<1$ it follows that there exists $0<q<1$ such that 
\Eq{*}{
  \limsup_{m\to\infty}\sqrt[m]{\|T^m\|^*_d}<q.
}
Then, there exists an $m_0\in\N$ such that, for all $m\geq m_0$,
\Eq{*}{
  \|T^m\|^*_d<q^m.
}
Therefore, for $n\geq m\geq m_0$, we have
\Eq{*}{
   \Bigg\|\sum_{k=0}^n T^k-\sum_{k=0}^m T^k\Bigg\|^*_d
   =\Bigg\|\sum_{k=m+1}^n T^k\Bigg\|^*_d
   \leq \sum_{k=m+1}^n\big\| T^k\big\|^*_d
   \leq \sum_{k=m+1}^nq^k\leq\frac{q^{m+1}}{1-q}.
}
This inequality yields that the sequence $S_n:=\sum_{k=0}^n T^k$ is a Cauchy sequence in $\E^d(X)$. Therefore it converges to an additive function $S$ belonging to $\E^d(X)$. One can easily see that $(I-T)\circ S_n=S_n\circ (I-T)=I-T^{n+1}$. Upon taking the limit $n\to\infty$, we get that $(I-T)\circ S=S\circ (I-T)=I$, which proves that $S$ is the inverse of $I-T$.
\end{proof}

\Cor{NIT}{Let $(X,+,d)$ be a complete metric Abelian group and let $T,S\in\E^d(X)$ such that $S$ is invertible with $S^{-1}\in\E^d(X)$ and $\min(\rho_d(T\circ S^{-1}),\rho_d(S^{-1}\circ T))<1$. Then $S-T$ is an invertible element of $\E^d(X)$.}

\begin{proof} If $\rho_d(T\circ S^{-1})<1$, then, by \thm{NIT}, $I-T\circ S^{-1}$ is an invertible element in $\E^d(X)$ and we have
\Eq{*}{
  S^{-1}\circ\big(I-T\circ S^{-1}\big)^{-1}
  = \big((I-T\circ S^{-1})\circ S\big)^{-1}
  =(S-T)^{-1}.
}
Therefore, $S-T$ is also an invertible element of $\E^d(X)$.

In the case when $\rho_d(S^{-1}\circ T)<1$ holds, then $I-S^{-1}\circ T$ is an invertible element of $\E^d(X)$ and we get
\Eq{*}{
  \big(I-S^{-1}\circ T\big)^{-1}\circ S^{-1}
  = \big(S\circ(I-S^{-1}\circ T)\big)^{-1}
  =(S-T)^{-1}.
}
Thus, again we obtain $S-T$ is an invertible element of $\E^d(X)$.
\end{proof}

\section{Main results}

Let $(X,+)$ be an Abelian group.
Given an endomorphism $T\in\E(X)$, we say that a subset $D\subseteq X$ is \emph{$T$-convex} if, for all $x,y\in D$, 
\Eq{*}{
  T(x)+(I-T)(y)\in D.
}
This condition is equivalent to the inclusion
\Eq{*}{
  T(D)+(I-T)(D)\subseteq D.
}
If $\T\subseteq\E(X)$, then a set $D\subseteq X$ is called \emph{$\T$-convex} if it is $T$-convex for all $T\in\T$. The class of $\T$-convex subsets of $X$ is denoted by $\C_\T(X)$ in what follows. In the particular case when $(X,+)$ is the additive group of a vector space and $T=tI$ for some $t\in[0,1]$, instead of $T$-convexity, we briefly speak about \emph{$t$-convexity} which is a commonly accepted notion (cf.\ \cite{Kuh84}). If $X$ is a uniquely $2$-divisible Abelian group, and $T=\frac12\Cdot I$, that is, $T(x):=\frac12\Cdot x$, then $T$-convex sets will also be termed \emph{midpoint convex}. One can immediately see that if the group $X$ is divisible by some $n\in \N$ and $T=\frac{1}{n}\Cdot I$, then $T$-convexity is equivalent to $n$-convexity defined in the previous section. 

The following lemma is obvious but useful.

\Lem{TC}{Given an endomorphism $T:X\to X$, a subset $D\subseteq X$ is $T$-convex if and only if, for all $p\in D$, 
\Eq{Tp}{
  T(D-p)\subseteq D-p.
}}

\begin{proof} The $T$ convexity of $D$ is equivalent to the property that, for all $p\in D$, the inclusion
\Eq{*}{
  T(D)+(I-T)(p)\subseteq D
}
holds which, by the additivity of $T$, is identical to the inclusion \eq{Tp}.
\end{proof}

The basic properties of $T$-convex sets are described in the following result.

\Thm{0}{For a nonempty subset $\T\subseteq\E(X)$, we have the following assertions.
\begin{enumerate}[(i)]
 \item $\C_\T(X)$ contains the empty set, the whole space $X$ and, for every $x\in X$, the singleton $\{x\}$.
 \item $\C_\T(X)$ is closed with respect to intersection and chain union.
 \item $\C_\T(X)$ is closed under algebraic addition.
 \item If $A\in\E(X)$ commutes with any member of $\T$ and $D\in\C_\T(X)$, then $A(D),A^{-1}(D)\in\C_\T(X)$.
\end{enumerate}}

\begin{proof} The first statement is obvious. 
The assertions concerning intersection and chain union are easy to verify. The closedness with respect to algebraic addition is a consequence of the additivity of the elements of $\T$. 

To prove the last assertion, let $A\in\E(X)$ commute with any member of $\T$ and let $D\in\C_\T(X)$.

First we prove that $A(D)$ is $\T$-convex. 
For this, let $u,v\in A(D)$. Then there exist $x,y\in D$ such that $A(x)=u$ and $A(y)=v$. Using that $D$ is $\T$-convex, we have that $T(x)+(I-T)(y)\in D$ for all $T\in\T$. Hence, for all $T\in\T$,
\Eq{*}{
  T(u)+(I-T)(v)=(T\circ A)(x)+((I-T)\circ A)(y)=A(T(x)+(I-T)(y))\in A(D),
}
which shows that $A(D)$ is $\T$-convex.

To prove that $A^{-1}(D)$ is $\T$-convex, let $x,y\in A^{-1}(D)$. Then $A(x),A(y)\in D$. By the $\T$-convexity of $D$, for all $T\in\T$, it follows that
\Eq{*}{
   A(T(x)+(I-T)(y))=T\circ A(x)+(I-T)\circ A(y)\in D.
}
Hence, $T(x)+(I-T)(y)\in A^{-1}(D)$ for all $T\in\T$. This completes the proof of the $\T$-convexity of $A^{-1}(D)$.
\end{proof}

In view of the closedness of $\C_\T(X)$ with respect to intersection, for every subset $S$ of $X$, the set 
\Eq{*}{
  \conv_\T(S):=\bigcap\{D\in\C_\T(X):S\subseteq D\}
}
is the smallest $\T$-convex set containing $S$, which will be called the \emph{$\T$-convex hull} of $S$.

Now, given a nonempty subset $D\subseteq X$, we consider the collection of endomorphisms $T$ of $X$ that make $D$ to be $T$-convex: 
\Eq{*}{
  \T_D:=\{T\in\E(X)\mid \mbox{$D$ is $T$-convex}\}.
}
It is obvious that, for every set $D$, we have $0,I\in\T_D$ and $0,I\in\T^d_D$ (if $X$ is a metric Abelian group). The next result describes a convexity property of $\T_D$.

\Thm{P1}{Let $D\subseteq X$ be a nonempty set. Then $\T_D$ is a $\widetilde\T_D$-convex subset of $\E(X)$. If $(X,+,d)$ is a metric Abelian group, then $\T^d_D$ is a $\widetilde\T^d_D$-convex subset of $\E^d(X)$. In particular, these sets are closed with respect to the composition of maps.}

\begin{proof} Let $T,T_1,T_2\in\T_D$ and set $S:=T\circ T_1+(I-T)\circ T_2$. Then, by the $T_1$-, $T_2$- and $T$-convexity of $D$, for all $x,y\in D$, we have
\Eq{*}{
  u:=T_1(x)+(I-T_1)(y)\in D, \qquad v:=T_2(x)+(I-T_2)(y)\in D,
}
hence 
\Eq{*}{
   T(u)+(I-T)(v)\in D.
}
On the other hand,
\Eq{*}{
  T(u)+(I-T)(v)
  &=T\big(T_1(x)+(I-T_1)(y)\big)+(I-T)\big(T_2(x)+(I-T_2)(y)\big)\\
  &=\big(T\circ T_1+(I-T)\circ T_2\big)(x)+(I-T\circ T_1-(I-T)\circ T_2)(y)
  =S(x)+(I-S)(y).
}
Therefore, $S\in\T_D$ follows. This yields that $\T_D$ is $\widetilde{T}$-convex for all $T\in\T_D$, which was to be proved.

The proof of the second assertion is completely analogous, the last statement follows by taking $T_2=0$ in the above argument. 
\end{proof}

\Exa{345}{
It is not true in general that an arbitrary set $\T\subseteq \E(X)$ is $\widetilde{\T}$-convex. To visualise this take $(X,+)=(\Z,+)$ and $\T=\{ \pi_3, \pi_4, \pi_5 \}$. Clearly $\pi_3\circ \pi_4 + (I-\pi_3)\circ \pi_5 =\pi_2 \notin \T$.
}

\Cor{1}{Let $D\subseteq X$ be a nonempty set. Then $\T_D$ and also $\T^d_D$ (if $(X,+,d)$ is a metric Abelian group) are closed under multiplication and under the mappings 
\Eq{maps}{
T\mapsto I-T\qquad\mbox{and}\qquad (T,S)\mapsto T\circ S+(I-T)\circ(I-S).
}}

\begin{proof} Let $T,T_1,T_2,S\in\T_D$. Then, by \thm{P1}, we have $T\circ T_1+(I-T)\circ T_2\in\T_D$. Putting $T_1:=S$, and $T_2:=0$, we get $T\circ S\in\T_D$. On the other hand, letting $T_1:=0$, $T_2:=I$, we obtain $I-T\in\T_D$. Finally, substituting $T_1:=S$ and $T_2:=I-S$, we can see that $T\circ S+(I-T)\circ(I-S)\in\T_D$ is also valid.

The proof for the metric Abelian group setting is completely analogous. 
\end{proof}

In the theorem below, we show that $T$-convexity implies midpoint convexity under certain conditions on $X$ and $T$. 


\Thm{2}{Let $(X,+,d)$ be a complete metric uniquely 2-divisible Abelian group and $T\in\E^d(X)$ such that $\rho_d(2\Cdot T-I)<1$. Then, for every nonempty $T$-convex set $D\subseteq X$, the set $\cl(\T^d_D)$ is a midpoint convex subset of $\E^d(X)$. Furthermore, every closed $T$-convex subset of $X$ is also midpoint convex.}

\begin{proof} Let $D$ be a nonempty $T$-convex subset of $X$. Define the sequence of operators $T_n$ by the following recursion:
\Eq{RTn}{
  T_1:=T,\qquad T_{n+1}:=T_n^2+(I-T_n)^2 \quad(n\in\N).
}
Then, in view of the last assertion of \cor{1}, we have that $T_n\in\T^d_D$ for all $n\in\N$. On the other hand it follows by induction that 
\Eq{Tn}{
  T_n=\frac{1}{2}\Cdot \Big(I+(2\Cdot T-I)^{2^{n-1}}\Big).
}
On can see that this formula is correct for $n=1$. Assume that it is valid for $n=k\geq 1$. Then
\Eq{*}{
  T_{k+1}=T_k^2+(I-T_k)^2
  &= \frac{1}{4}\Cdot\Big(I+(2\Cdot T-I)^{2^{k-1}}\Big)^2
    +\frac{1}{4}\Cdot\Big(I-(2\Cdot T-I)^{2^{k-1}}\Big)^2 \\
  &=\frac14\Cdot\Big(I+2\Cdot(2\Cdot T-I)^{2^{k-1}}+(2\Cdot T-I)^{2^{k}}
      +I-2\Cdot(2\Cdot T-I)^{2^{k-1}}+(2\Cdot T-I)^{2^{k}}\Big)\\
  &=\frac{1}{2}\Cdot \Big(I+(2\Cdot T-I)^{2^{k}}\Big),  
}
which proves the validity of \eq{Tn} for $n=k+1$.

From the assumption $\rho_d(2\Cdot T-I)<1$, using \lem{SR} we obtain that $T_n$ converges to $\frac12\Cdot I$. If $R,S\in\cl\T^d_D$, then there exist sequences $R_n$ and $S_n$ in $\T^d_D$ converging to $R$ and $S$, respectively. By \thm{P1}, we have that $T_n\circ R_n+(I-T_n)\circ S_n\in\T^d_D$. Upon taking the limit, it follows that $\frac12\Cdot (R+S)\in\cl\T^d_D$ proving that the set $\cl\T^d_D$ is midpoint convex.

To complete the proof, assume that $D$ is a closed $T$-convex subset of $X$. Then it is also $T_n$-convex for all $n\in\N$. Then, taking the limit $n\to\infty$, it follows that $D$ is midpoint convex.
\end{proof}

\Thm{nk}{Let $(X,+,d)$ be a metric Abelian group, assume that there exists $n_0\in\N$ such that $\mu_d(n_0)>1$. Assume that either $X$ is complete or $n_0\Cdot X$ is closed and let $D$ be a bounded $n_0$-convex set. Then, for all $n\in\N$ and $T_1,\dots,T_n\in\T^d_D$, 
\Eq{*}{
   T_1(D)+\cdots+T_n(D)\subseteq \cl\big((T_1+\cdots+T_n)(D)\big).
}}

\begin{proof} If $D$ is $n_0$-convex, then it is also $n_0^k$-convex for all $k\in\N$. Let $n\in\N$, $T_1,\dots,T_n\in\T^d_D$ and then choose $k\in\N$ so that $n_0^k>n$.

By using the $n_0^k$- and $T_1$-,\dots, $T_n$-convexity of $D$,
we arrive at an inclusion where \thm{RCT} can be applied: 
\Eq{*}{
  \sum_{i=1}^n T_i(D)&+\bigg(\sum_{i=1}^n (I-T_i)+(n_0^k-n)I\bigg)(D)
  \subseteq \sum_{i=1}^n \big(T_i(D)+(I-T_i)(D)\big)+[n_0^k-n](D) \\
  &= \sum_{i=1}^n D+[n_0^k-n](D)=[n_0^k]D=n_0^k\Cdot D
  =\bigg(\sum_{i=1}^n T_i+\sum_{i=1}^n (I-T_i)+(n_0^k-n)I\bigg)(D)\\
  &\subseteq\bigg(\sum_{i=1}^n T_i\bigg)(D)+\bigg(\sum_{i=1}^n (I-T_i)+(n_0^k-n)I\bigg)(D).
}
Now, taking
\Eq{ABC}{
  A:=\sum_{i=1}^n T_i(D),\qquad B:=\bigg(\sum_{i=1}^n T_i\bigg)(D)
  \qquad\mbox{and}\qquad C:=\bigg(\sum_{i=1}^n (I-T_i)+(n_0^k-n)I\bigg)(D),
}
we can see that $A+C\subseteq B+C\subseteq\cl(B)+C$ holds. Furthermore, due to the Lipschitz property of the endomorphisms $T_1,\dots,T_n$ and the boundedness of $D$, the set $C$ is bounded. The $n_0$-convexity of $D$ easily implies that $B$ is $n_0$-convex. Using that $\mu_d(n_0)>0$, \lem{nx} implies 
\Eq{*}{
  [n_0]\cl(B)\subseteq \cl\big([n_0]B\big)
  \subseteq \cl\big(n_0\Cdot B\big)
  \subseteq \cl\big(n_0\Cdot \cl(B)\big)
  = n_0\Cdot \cl(B).
}
Thus, $\cl(B)$ is also $n_0$-convex. Hence, \thm{RCT} implies that $A\subseteq \cl(B)$, which was to be proved. 
\end{proof}

With additional properties on the set $D$, or on the space $X$, or on the operator $T_1+\dots+T_n$, we have the following stronger statement.

\Thm{nk+}{Let $(X,+,d)$ be a metric  Abelian group, let $n_0\in\N$ such that $\mu_d(n_0)>1$ and let $D$ be a closed bounded $n_0$-convex set. Let $n\in\N$ and let $T_1,\dots,T_n\in\T^d_D$. Assume that at least one of the following conditions holds:
\begin{enumerate}[(i)]
 \item $D$ is compact;
 \item $X$ is complete and $\mu_d(T_1+\dots+T_n)>0$;
 \item $(T_1+\dots+T_n)(X)$ is closed and $\mu_d(T_1+\dots+T_n)>0$.
\end{enumerate}
Then 
\Eq{TT}{
   T_1(D)+\cdots+T_n(D)\subseteq (T_1+\cdots+T_n)(D).
}}

\begin{proof} As we have seen it in the proof of \thm{nk}, with notation \eq{ABC}, we have $A+C\subseteq B+C$. The set $B$ is $n_0$-convex. In order to deduce the inclusion $A\subseteq B$ with the aid of \thm{RCT}, it suffices to show that in each of the cases (i), (ii), (iii), the set $B$ is closed.

If assumption (i) is valid, then the continuity of the endomorphisms, implies that $B$ is compact and hence it is a closed set. 

If either (ii) or (iii) holds, then \lem{nx} with operator $T:=T_1+\dots+T_n$ implies that $B=T(D)$ is closed.
\end{proof}

\Cor{nkc1}{Let $(X,+,d)$ be a metric group, let $n_0\in\N$ such that $\mu_d(n_0)>1$ and let $D$ be a compact $n_0$-convex set. Then $D$ is $n$-convex for all $n\in\N$.}

\begin{proof} The statement immediately follows by applying the first assertion of \thm{nk+} for the setting $T_1=\dots=T_n=I$. 
\end{proof}

The above theorem also implies a further closure property of the set $\T_D^d$. 

\Cor{nkc2}{Let $(X,+,d)$ be a metric  Abelian group, let $n_0\in\N$ such that $\mu_d(n_0)>1$ and let $D$ be a closed bounded $n_0$-convex set. Let $n\in\N$ and $T_1,\dots,T_n\in\T^d_D$ be such that $T:=T_1+\cdots+T_n$ is a bijection with $T^{-1}\in\E^d(X)$. Then, for all $k\in\{1,\dots,n-1\}$, we have $T^{-1}\circ(T_1+\cdots+T_k)\in\T_D^d$.}

\begin{proof} Provided that $T$ is a bijection with $d$-bounded inverse, we have that assumption (iii) of \thm{nk+} is satisfied, thus inclusion in \eq{TT} is valid. This implies that
\Eq{*}{
   (T^{-1}\circ T_1)(D)+\cdots+(T^{-1}\circ T_n)(D)\subseteq D.
}
Hence, for all $k\in\{1,\dots,n-1\}$,
\Eq{*}{
  \bigg(\sum_{i=1}^k T^{-1}\circ T_i\bigg)(D)
  +\bigg(I-\sum_{i=1}^k T^{-1}\circ T_i\bigg)(D)
  =\bigg(\sum_{i=1}^k T^{-1}\circ T_i\bigg)(D)
  +\bigg(\sum_{i=k+1}^n T^{-1}\circ T_i\bigg)(D)\subseteq D,
}
which shows that $T^{-1}\circ(T_1+\cdots+T_k)\in\T_D^d$. 
\end{proof}


\end{document}